\documentclass[10pt,twoside]{article}

\usepackage[centertags]{amsmath}
\usepackage{amsfonts}
\usepackage{amssymb}
\usepackage{amsthm}
\usepackage{epsfig}
\usepackage{color}
\allowdisplaybreaks[4]
\date{}

\begin{document}
\baselineskip 15pt \setcounter{page}{1}
\title{\bf \Large  The asymptotic distribution of maxima of stationary random sequences under random replacing
\thanks{Research supported by Innovation of Jiaxing City: a program to support the talented persons
and Project of new economy research center of Jiaxing City (No. WYZB202254).}}
\author{{\small Yuwei Li$^{1,2}$, Zhongquan Tan$^2$\footnote{ E-mail address:  tzq728@163.com }}\\
\\
{\small\it
1. Department of Mathematic, Zhejiang Normal University, Jinhua 321004, PR China}\\
{\small\it
2. College of Data Science, Jiaxing University, Jiaxing 314001, PR China}\\
}
 \maketitle
 \baselineskip 15pt

\begin{quote}
{\bf Abstract:}\ \ In this paper, we investigated  the
effect on extreme of random replacing for a stationary sequence satisfying a type of long dependent condition and a local dependent condition,
and derived the joint asymptotic distribution of maximum from the stationary sequence and the maximum from the random replacing sequence.
We also provided several applications for our main results.

{\bf Key Words:}\ \ extreme value theory; random replacing; asymptotic distribution; stationary sequences

{\bf AMS Classification:}\ \ Primary 60G70; secondary 60G10

\end{quote}

\section{Introduction}

In practice, missing data may occur randomly. What can we do when a missing data occur? In the majority of situations, we have two ways to deal with it:
regarding the data sample as a smaller sample with random sample size or replacing the missing data by a fixed value.
When studying the properties of the underlying random variables it may be crucial to take into account
the effect of the missing data in the former case and the replacing data in the latter case.
In extreme value theory, these effects have been extensively studied in the literature.

Suppose that $\{X_n,n\geq 1\}$ is a sequence of stationary random variables with the marginal distribution function $F(x)$ and some of the random variables in the sequence are missing randomly. Let $\boldsymbol{\varepsilon}= \{\varepsilon_n,n\geq 1\}$ be a sequence of Bernoulli random variables which captures the phenomenon of random missing.

For the random sequence $\{X_n,n\geq 1\}$, define its random missing sequence as:
\begin{eqnarray}
\label{M1}
\widetilde{X_{n}}(\boldsymbol{\varepsilon})=\varepsilon _nX_n+(1-\varepsilon_n)x_{F},
\end{eqnarray}
 where $x_{F}=\inf\{x\in \mathbb{R}: F(x)>0\}$.
Suppose that $\{\varepsilon_n,n\geq 1\}$ is independent of $\{X_{n},{n\geq 1}\}$ and let $S_{n}=\sum_{k\leq n}\varepsilon_{k}$ be the numbers of the observed variables satisfying
\begin{eqnarray*}
\frac{S_{n}}{n}\stackrel{P}{\longrightarrow} \lambda,\ \ \mbox{as}\quad n\rightarrow\infty,
\end{eqnarray*}
where $\lambda$ is a random or nonrandom variable.

When $\lambda\in [0,1]$ is a constant, Mladenovi\'{c} and Piterbarg (2006) first studied the effect of the missing data on extremes under conditions $D(u_{n},v_{n})$ and $D'(u_{n})$ (see the definitions in Section 2), and derived the asymptotic distribution of maximum from stationary sequences and the maximum from the random missing sequences. More precisely, for any $x < y \in \mathbb{R}$, they proved
\begin{eqnarray}
\label{MP}
\lim_{n\rightarrow\infty}P\left(M_{n}(\widetilde{X}(\boldsymbol{\varepsilon}))\leq a_{n}^{-1}x+b_{n}, M_{n}(X)\leq a_{n}^{-1}y+b_{n}\right)=G^{\lambda}(x)G^{1-\lambda}(y),
\end{eqnarray}
with $a_{n}>0$ and $b_{n}\in \mathbb{R}$, where $M_{n}(\widetilde{X}(\boldsymbol{\varepsilon}))=\max\{\widetilde{X_{k}}(\boldsymbol{\varepsilon}), k=1,2,\ldots,n\}$ and $M_{n}(X)=\max\{X_{k}, k=1,2,\ldots,n\}$.
For the extensions of (\ref{MP}), we refer to Cao and Peng (2011) and Peng et al. (2010) for the Gaussian cases; Glava\u{s} et al. (2017) for autoregressive process; Glava\u{s} and Mladenovi\'{c} (2020) for linear process; Panga and Pereira (2018) for nonstationary random fields; and Peng et al. (2010), Tong and Peng (2011) and Tan and Wang (2012) for the almost sure limit theorem.

When $\lambda\in [0,1]$ a.s.  is  a random variable, Krajka (2011) obtained the following result: for any $x< y \in \mathbb{R}$
\begin{eqnarray}
\label{KR}
\lim_{n\rightarrow\infty}P\left(M_{n}(\widetilde{X}(\boldsymbol{\varepsilon}))\leq a_{n}^{-1}x+b_{n}, M_{n}(X)\leq a_{n}^{-1}y+b_{n}\right)= E[G^{\lambda}(x)G^{1-\lambda}(y)].
\end{eqnarray}
Hashorva et al. (2013) extended the results of (\ref{KR}) to weakly and strongly dependent Gaussian sequences. For more related studies under this situation, we refer to  Krajka and Rychlik (2014) and Zheng and Tan (2022).

For the random sequence $\{X_n,n\geq 1\}$, define its random replacing sequence as:
\begin{eqnarray}
\label{M2}
X_{n}(\boldsymbol{\varepsilon})=\varepsilon _nX_n+(1-\varepsilon_n)\widehat{X_n},
\end{eqnarray}
where the sequence $\{\widehat{X_n}, n\geq 1\}$ is an independent copy of $\{X_n, n\geq 1\}$.
Under a type of distribution mixing condition,  Hall and H\"{u}sler (2006) first investigated the effect on extreme index of random replacing data, when $\{\varepsilon_n,n\geq 1\}$ is independent of $\{X_{n},{n\geq 1}\}$. To our knowledge,  the limiting properties describing the asymptotic relation between the extremes of some sequences and their replaced ones have not been well studied in literature before.

The objective of this paper is to continue to study the model (\ref{M2}) and  to derive the joint asymptotic distribution of maximum from some stationary sequences and the maximum from its random replacing sequences.
For the strongly mixing random sequences, when $\{\varepsilon_n,n\geq 1\}$ is dependent on $\{X_{n},{n\geq 1}\}$, this problem was considered by Robert (2010).
Robert (2010) proved the maximum from the strongly mixing random sequences and the maximum from its random replacing sequences are asymptotically dependent.
The strongly mixing condition is so restrictive in extreme value theory that it excludes the weakly dependent Gaussian sequences.
In fact, we will show under very mild conditions (conditions $D(u_{n},v_{n})$ and $D'(u_{n})$) that the maximum from a stationary random sequences and the maximum from its random replacing sequences are asymptotic dependent. Several applications on  Gaussian sequences, $\chi$ sequences and Gaussian order statistics sequences are also given to illustrate  the results.

The rest of the paper is organized as follows.  The main results of the paper and some applications are given in Section 2 and their proofs are collected in Section 3.

\section{Main results}

We first give the long dependent condition $D(u_{n},v_{n})$ and the local dependent condition $D'(u_{n})$.

\textbf{Definition 2.1}. Let $\{X_n, n\geq 1\}$ be a strictly stationary random sequence and $u_{n}$ and $v_{n}$ be two sequences of real numbers.
The condition $D(u_n, v_n)$ is satisfied, if for all $A_1, A_2, B_1, B_2 \subset \{1,2, \dots, n\}$, such that
$\max_{a\in A_1\cup A_2, b\in B_1\cup B_2 } \mid b-a\mid\ge l, A_1\cap  A_2=\phi,  B_1\cap B_2=\phi$, the following inequlaity holds:
\begin{eqnarray*}
&&\bigg|P\left[\bigcap_{j\in A_1\cup B_1}^{} \{X_j\le u_{n}\}\cap \bigcap_{j\in A_2\cup B_2}\{X_j\le v_{n}\}\right]
\\&&\ \ \ -P\left[\bigcap_{j\in A_1}\{X_j\le u_n \}\cap \bigcap_{j\in A_2}\{X_j\le v_n \}\right]P\left[\bigcap_{j\in B_1}^{} \{X_j\le u_n \}\cap \bigcap_{j\in B_2}\{X_j\le v_n \}\right]\bigg|
\\&& \leq \alpha_{n,l},
\end{eqnarray*}
and $\alpha_{n,l_n} \rightarrow 0$ as  $n  \rightarrow\infty$  for some $l_n =o(n)$.

Condition $D(u_{n},v_{n})$ is taken from Mladenovi\'{c} and Piterbarg (2006), which is an extension of the classical long dependent condition $D(u_{n})$ in Leadbetter et al. (1983).
It is easy to check that independent and identically distributed random sequence, $m$-dependent random sequence and strongly dependent random sequence satisfy
Condition $D(u_n, v_n)$.

The following local dependent condition $D'(u_{n})$ taken from Leadbetter et al. (1983) is an anti-cluster condition, which bounds the probability of more than one exceedance above the levels $u_{n}$ in a interval with a few indexes.

\textbf{Definition 2.2}. Let $\{X_n, n\geq 1\}$ be a strictly stationary random sequence of random variables and $\{u_n, n\geq 1\}$  be a sequence of real numbers. We say that $\{X_n, n\geq 1\}$ satisfy the Condition $D'(u_n)$ if
\begin{equation*}
\lim_{k \to \infty} \limsup_{n \to \infty} n\sum_{j=2}^{\lfloor n/k\rfloor}  P[X_1>u_n, X_j>u_n]=0,
\end{equation*}
where $\lfloor x\rfloor $ denotes the integral part of $x$.
Obviously, if $\{X_n, n\geq 1\}$ is a sequence of independent and identically distributed random variables with $\lim_{n \to \infty} n P[X_1>u_n]=c$, then Condition $D'(u_n)$ holds.

\textbf{Theorem 2.1}. {\sl Suppose that the following conditions are satisfied:

(a) $\{X_n, n\geq 1\}$ is a strictly stationary random sequence with the marginal distribution function $F(x)$ satisfying conditions $D(u_n, v_n)$  and  $D'(u_n)$  for
 $u_n= a_n^{-1}x + b_n$ and $v_n= a_n^{-1}y + b_n$;

(b) $F \in D(G)$, i.e., for some real constants $a_n > 0 , b_n\in \mathbb{R}$, and any $x\in \mathbb{R}$ the equality
\begin{eqnarray}
\label{2.10}
\lim_{n \to \infty} F^n(a_n^{-1}x+b_n)=G(x),
\end{eqnarray}holds for every continuity point of G;

(c) $\boldsymbol{\varepsilon}= \{\varepsilon_n, n\geq 1\}$ is a sequence of Bernoulii random variables and is independent of $\{X_n, n\geq 1\}$ such that
 \begin{eqnarray}
 \label{2.1}
\frac{{S_n }}{n} \stackrel{p}\longrightarrow \lambda,\ \ \mbox{ as}\ \  n\rightarrow\infty,
 \end{eqnarray}
for some random variable $\lambda\in [0,1]$ a.s. Then  for any $x,y\in \mathbb{R}$
\begin{eqnarray}
 \label{T2.1}
 \lim_{n \to \infty} P[M_n(X(\boldsymbol{\varepsilon} ))\le a_n^{-1}x + b_n, M_n(X)\le a_n^{-1}y + b_n]=G(\min\{x, y\})EG^{1-\lambda }(\max\{x, y\}).
 \end{eqnarray} }

\textbf{Corollary 2.1}. {\sl Suppose the conditions of Theorem 2.1 hold. We have for any $x\in \mathbb{R}$
\begin{eqnarray}
 \label{C2.1}
 \lim_{n \to \infty} P[M_n(X(\boldsymbol{\varepsilon} ))\le a_n^{-1}x + b_n]=G(x)
 \end{eqnarray}
 and
 \begin{eqnarray}
 \label{C2.2}
 \lim_{n \to \infty} P[M_n(X)\le a_n^{-1}x+ b_n]=G(x).
 \end{eqnarray}}

\textbf{Remark 2.1}. Theorem 2.1 shows that the maxima of a stationary sequence and the maxima of its random replacing sequence are asymptotic
dependent. Corollary 2.1 shows that the random replacing does not affect the asymptotic distribution of the maxima.

\textbf{Remark 2.2}. (i). In model (\ref{M2}), the sequence $\{\widehat{X_n}, n\geq 1\}$ is an independent copy of $\{X_n, n\geq 1\}$. In some situations,
we may use an independent and identically distributed random sequence to replace the missing data. Since the dependence assumed  on the sequence in this paper is very weak, it is not hard to see that the results of Theorem 2.1 still hold in this case. \\
(ii). It is very interesting to discuss the case that the original random sequence and its random replacing sequence have different dependent structures.
It maybe lead to very different asymptotic distribution.  This will be studied for stationary Gaussian case in a forth coming paper.

Next, we give several applications of Theorem 2.1.

\textbf{Theorem 2.2}. {\sl Suppose that the following conditions are satisfied:

(a) $\{Y_n, n\geq 1\}$ is a sequence of standard (with 0 mean and 1 variance) Gaussian random variables  with covariance function $r_{n}=Cov(Y_{1}, Y_{n+1})$.
Assume that $r_{n}$ satisfies $r_{n}\log n\rightarrow0$ as $n\rightarrow\infty$.

(b) $\boldsymbol{\varepsilon}= \{\varepsilon_n, n\geq 1\}$ is a sequence of Bernoulii random variables and is independent of $\{Y_n, n\geq 1\}$ such that
 (\ref{2.1}) hols for some random variable $\lambda\in [0,1]$ a.s. Then for any $x,y\in \mathbb{R}$
\begin{eqnarray}
 \label{2.2}
&& \lim_{n \to \infty} P[M_n(Y(\boldsymbol{\varepsilon} ))\le a_n^{-1}x + b_n, M_n(Y)\le a_n^{-1}y + b_n]\nonumber\\
&&\ \ \ \ \ \ \ \ \ \ \ \ \ \ =\exp\left(-e^{-\min\{x, y\}}\right)E\exp\left(-(1-\lambda )e^{-\max\{x, y\}}\right),
 \end{eqnarray}
where $a_{n}=(2\log n)^{1/2}$ and $b_{n}=a_{n}-\frac{\log\log n+\log (4\pi)}{2a_{n}}$.
}

Let $\{Y_{n}, n\geq1\}$  defined as in Theorem 2.2 and let
$\{Y_{nj}, n\geq1\}, j=1,2,\cdots,d$ be its independent copies with $d\geq 1$.
Define two type of Gaussian functions as:
$$\chi_{n}=(\sum_{j=1}^{d}Y_{nj}^{2})^{1/2},\ \ n\geq1$$
and for $r\in\{1,2,\ldots,d\}$
$$O_{n}^{(d)}:=\min_{j=1}^{d}Y_{nj}\leq \cdots\leq O_{n}^{(r)}\leq \cdots\leq O_{n}^{(1)}:=\max_{j=1}^{d}Y_{nj},\ \ n\geq1.$$
Note that $\{\chi_{n}, n\geq1\}$ is a sequence of chi-random variables and
$\{O_{n}^{(r)}, n\geq1\}$ is called Gaussian $r$-th order statistics random variables.
The above two type of Gaussian functions play a crucial role in various statistical applications. For more details, we refer to
Tan and Hashorva (2013a, 2013b) for chi-random variables and
D\c{e}bicki et al. (2015, 2017) and Tan (2018) for Gaussian $r$-th order statistics.

\textbf{Theorem 2.3}. {\sl Suppose the conditions of Theorem 2.2 hold. We have for any $x,y\in \mathbb{R}$
\begin{eqnarray}
 \label{2.3}
&& \lim_{n \to \infty} P[M_n(\chi(\boldsymbol{\varepsilon} ))\le a_n^{-1}x + b_n, M_n(\chi)\le a_n^{-1}y + b_n]\nonumber\\
&&\ \ \ \ \ \ \ \ \ \ \ \ \ \ =\exp\left(-e^{-\min\{x, y\}}\right)E\exp\left(-(1-\lambda )e^{-\max\{x, y\}}\right),
 \end{eqnarray}
where
$a_{n}=(2\log n)^{1/2}$ and $b_n=a_n+ a_n^{-1}\log \Bigl( 2^{1- d/2}(\Gamma(d/2))^{-1}a_n^{d-2}\Bigr)$
with $\Gamma(\cdot)$ the Euler gamma function.
}

\textbf{Theorem 2.4}. {\sl Suppose the conditions of Theorem 2.2 hold. We have for any $x,y\in \mathbb{R}$
\begin{eqnarray}
 \label{2.4}
&& \lim_{n \to \infty} P[M_n(O^{(r)}(\boldsymbol{\varepsilon} ))\le a_n^{-1}x + b_n, M_n(O^{(r)})\le a_n^{-1}y + b_n]\nonumber\\
&&\ \ \ \ \ \ \ \ \ \ \ \ \ \ =\exp\left(-e^{-\min\{x, y\}}\right)E\exp\left(-(1-\lambda )e^{-\max\{x, y\}}\right),
 \end{eqnarray}
where $a_{n}=(2r\log n)^{1/2}$ and $b_{n}=\frac{1}{r}\alpha_{n}+ \alpha_{n}^{-1}\log \Bigl(\alpha_{n}^{-r} C_{d}^{r}(2\pi)^{-r/2}\Bigr)$ with $C_{d}^{r}=\frac{d!}{r!(d-r)!}$.
}

\section{Proofs}
We first state and prove several lemmas which will help us to complete the proof of the Theorem 2.1.
For the simplicity, let in the following part $u_n(x)= a_n^{-1}x + b_n$ and $u_n(y)= a_n^{-1}y + b_n$.
Let $ \boldsymbol{\alpha}=\{\alpha_n, n\geq 1\}$ be a sequence of 0 and 1 ($\boldsymbol{\alpha} \in\{0, 1\}^\mathbb{N}$).
For the arbitrary random or nonrandom sequence $\boldsymbol{\beta} =\{\beta_n,n\geq 1\}$ of 0 and 1 and subset $I\subset \mathbb{N}$, put
$$M(X(\boldsymbol{\beta} ),I)=\max \{X_i(\boldsymbol{\beta}),  i\in I\},\ \ M(X,I)=\max \{X_i,  i\in I\}.$$

\textbf{Lemma 3.1}. {\sl Let $I_1, I_2, \cdots , I_k$ be subsets of \{1,2, \dots, n\} such that $\min I_m - \max I_s \geq l $, for $ k \geq m >s \geq 1$. Under the conditions of Theorem 2.1, we have
\begin{eqnarray*}
&&\bigg| P\left[\bigcap_{s=1}^{k} M(X(\boldsymbol{\alpha}  ), I_s)\le u_n(x), M(X, I_s)\le u_n(y)\right]-
 \prod_{s=1}^{k} P\left[M(X(\boldsymbol{\alpha}  ), I_s) \le u_n(x), M(X, I_s) \le u_n(y)\right]\bigg| \\
&&\leq   2(k-1) \alpha_{n,l},
\end{eqnarray*}
where $\alpha_{n,l}\rightarrow 0 $ as $n \rightarrow \infty$, for some $l_n =o(n)$.
}

\textbf{Proof}: Define
 $I_s^* = \{i:  i\in I_s, \alpha_i=1\},$
$I_s^{**} = \{i:  i\in I_s, \alpha_i=0\}$, obviously $I_s=I_s^*\cup I_s^{**}.$
By the independence between $\{\widehat{X_i},i\ge 1\}$ and $\{X_i, i\ge 1\}$, we have
\begin{eqnarray*}
&&\bigg| P\left[\bigcap_{s=1}^{k} M(X(\boldsymbol{\alpha}  ), I_s)\le u_n(x), M(X, I_s)\le u_n(y)\right]- \prod_{s=1}^{k} P\left[M(X(\boldsymbol{\alpha}  ), I_s) \le u_n(x), M(X, I_s) \le u_n(y)\right]\bigg|  \\
&&=  \bigg|P\left[\bigcap_{s=1}^{k}\{ M(\widehat{X},I_s^{**}) \le u_n(x)\}\right]P\left[\bigcap_{s=1}^{k}\{M(X,I_s^{*}) \le u_n(x), M(X,{I_s}) \le u_n(y)\}\right] \\
 &&\ \ \ \ - \prod_{s=1}^{k} P\left[M(\widehat{X},I_s^{**}) \le u_n(x)\right]P\left[ M(X,I_s^{*}) \le u_n(x), M(X,{I_s})\le u_n(y)\right]\bigg| \\
 &&\leq  \bigg|P\left[\bigcap_{s=1}^{k}\{ M(X,I_s^{*})\le u_n(x), M(X,{I_s})\le u_n(y)\}\right]-\prod_{s=1}^{k} P\left[M(X,I_s^{*})\le u_n(x) , M(X,{I_s})\le u_n(y)\right]\bigg|  \\
  &&\ \ \ \ + \bigg|P\left[ \bigcap_{s=1}^{k}\{ M(\widehat{X},I_s^{**})\le u_n(x) \}\right]-\prod_{s=1}^{k} P\left[M(\widehat{X},I_s^{**}) \le u_n(x)\right]\bigg|.
\end{eqnarray*}
Using the method of mathematical induction, by the similar arguments as the proof of Lemma 4.2 in Mladenovi\'{c} and Piterbarg (2006),  we get
\begin{eqnarray*}
 &&\bigg|P\left[\bigcap_{s=1}^{k}\{ M(X,I_s^{*})\le u_n(x), M(X,{I_s})\le u_n(y)\}\right]-\prod_{s=1}^{k} P\left[M(X,I_s^{*})\le u_n(x) , M(X,{I_s})\le u_n(y)\right]\bigg| \\
 &&\leq (k-1) \alpha_{n,l},
\end{eqnarray*}
where $\alpha_{n,l}\rightarrow 0 $ as $n \rightarrow \infty$, for some $l_n =o(n)$.
Note that $D(u_n, v_n)$ implies $D(u_n)$ by letting $v_n\rightarrow\infty$.
Similarly, using Condition $D(u_n)$, we get
\begin{eqnarray*}
\bigg|P\left[ \bigcap_{s=1}^{k}\{ M(\widehat{X},I_s^{**})\le u_n(x) \}\right]-\prod_{s=1}^{k} P\left[M(\widehat{X},I_s^{**}) \le u_n(x)\right]\bigg|
\leq (k-1) \alpha_{n,l},
\end{eqnarray*}
where $\alpha_{n,l}\rightarrow 0 $ as $n \rightarrow \infty$, for some $l_n =o(n)$. The proof of the lemma is complete.

For any positive integer $n$, let us define $\mathbf{N_n}=\left \{1, 2, \cdots, n \right \}$ and let $k$ be a fixed positive integer and $t=\lfloor \frac{n}{k}\rfloor $.
Let $$K_s=\{i: (s-1)t+1\leq i \leq st\}, \ \ s=1,2,\ldots, k.$$

\textbf{Lemma 3.2}.  {\sl Under the conditions of Theorem 2.1, we have
\begin{eqnarray}
\label{Lem3.2.1}
&&\bigg|P[M_n(X(\boldsymbol{\alpha} ))\le u_n(x), M_n(X)\le u_n(y)]-
 \prod_{s=1}^{k} P[M(X(\boldsymbol{\alpha}  ), K_s)\le u_n(x), M(X, K_s)\le u_n(y)]\bigg|\nonumber \\
&&\leq  2(k-1) \alpha_{n,l} + (\frac{4kl}{n}+\frac{2}{k})  n [1-F(u_n(\min\{x,y\}))],
\end{eqnarray}
where $\alpha_{n,l}\rightarrow 0 $ as $n \rightarrow \infty$, for some $l_n =o(n)$.
}

\textbf{Proof}: Divide each $K_s$ into two subsets
$I_s=\{i: (s-1)t+1\leq i \leq st-l\}$ and $J_s=\{i: st-l+1\leq i \leq st\}$, $s=1,2,\ldots, k$ and
denote also $J_{k+1}=\mathbf{N_n}-\cup_{s=1}^{k}K_s$. Obviously,  the absolute in  (\ref{Lem3.2.1}) is bounded above by
\begin{eqnarray*}
 &&\left|P\bigg[\bigcap_{s=1}^{k} M(X(\boldsymbol{\alpha} ), I_s)\le u_n(x), M(X, I_s) \le u_n(y)]\bigg]-P\bigg[M_n(X(\boldsymbol{\alpha} ))\le u_n(x), M_n(X)\le u_n(y)\bigg]\right|\\
  &&+\left|P\bigg[\bigcap_{s=1}^{k}\{ M(X(\boldsymbol{\alpha}), I_s)\le u_n(x), M(X, I_s) \le u_n(y)\}\bigg]-\prod_{s=1}^{k} P\bigg[M(X(\boldsymbol{\alpha}   ), I_s) \le u_n(x),  M(X, I_s) \le u_n(y)\bigg]\right|\\
  &&+\left|\prod_{s=1}^{k} P\bigg[M(X(\boldsymbol{\alpha}   ), I_s) \le u_n(x),  M(X, I_s) \le u_n(y)\bigg]-\prod_{s=1}^{k} P\bigg[M(X(\boldsymbol{\alpha}   ), K_s) \le u_n(x),  M(X, K_s) \le u_n(y)\bigg]\right|\\
 &&=P_{n1}+P_{n2}+P_{n3}.
\end{eqnarray*}
Note that $\sharp(J_{k+1})<t$, where $\sharp(A)$ denotes the cardinality of the set $A$. By the stationarity of $\{X_i, i\ge 1\}$ and $\{\widehat{X_i}, i\ge 1\}$, we  get
\begin{eqnarray*}
P_{n1}&\leq& kP[ M(X(\boldsymbol{\alpha} ), J_1)>u_n(x)]+P[ M(X(\boldsymbol{\alpha} ), J_{k+1})>u_n(x)]\\
&+&kP[M(X, J_1)>u_n(y)]+ P[M(X, J_{k+1})>u_n(y)]\\
&\leq& kl(1-F(u_n(x)))+t(1-F(u_n(x)))+ kl(1-F(u_n(y)))+t(1-F(u_n(y)))\\
&\leq& 2kl(1-F(u_n(\min\{x,y\})))+2t(1-F(u_n(\min\{x,y\}))).
\end{eqnarray*}
Similarly, we have
\begin{eqnarray*}
P_{n3}&\leq& 2kl(1-F(u_n(\min\{x,y\}))).
\end{eqnarray*}
From Lemma 3.1, we know
$$P_{n2} \leq 2(k-1) \alpha_{n,l}.$$
The proof is complete.

Let
$$A_{sj}=\{X_{(s-1)t+j}> u_n(\min\{x,y\})\},\ \ j=1,2,\ldots, t.$$

\textbf{Lemma 3.3}. {\sl Under the conditions of Theorem 2.1, for any $0 \leq r \leq 2^k-1$, we have
\begin{eqnarray*}
&&1-t[1-F(u_n(\min\{x,y\}))]-(1-\frac{r}{2^k} )t[1-F(u_n(\max\{x,y\}))]\\
&&\ \ \ +\left  [  \frac{\sum_{j\in K_s}^{}\alpha_j  }{t} -\frac{r}{2^k}\right ]t[1-F(u_n(\max\{x,y\}))]\\
&&\leq P[M(X(\boldsymbol{\alpha}  ), K_s)\le u_n(x), M(X, K_s)\le u_n(y)]\\
&&\leq 1-t[1-F(u_n(\min\{x,y\}))]-(1-\frac{r}{2^k} )t[1-F(u_n(\max\{x,y\}))]\\
&&\ \ \ +\left  [  \frac{\sum_{j\in K_s}^{}\alpha_j  }{t} -\frac{r}{2^k}\right ]t[1-F(u_n(\max\{x,y\}))]\\
&&\ \ \ +2t(1-F(u_n(\min\{x,y\})))^2+t\sum_{j=2}^{t} P[A_{s1}, A_{sj}].
\end{eqnarray*}
}

\textbf{Proof}:
Define  $K_s^* = \{i: i\in K_s, \alpha_i=1\}$ and
$K_s^{**} = \{i: i\in K_s, \alpha_i=0\}$, so $K_s=K_s^*\cup K_s^{**}$, $ s=1,2,\cdots, k.$
By the Bonferoni inequality and the independence between $\{\widehat{X_i},i\ge 1\}$ and $\{X_i, i\ge 1\}$,  we have
\begin{eqnarray*}
&&P[M(X(\boldsymbol{\alpha}), K_s)\le u_n(x), M(X, K_s)\le u_n(y)]\\
&&=P[M(\widehat{X}, K_s^{**})\le u_n(x), M(X, K_s^{*})\le u_n(\min\{x,y\}), M(X, K_s^{**})\le u_n(y)]\\
&&\leq 1-P[M(\widehat{X}, K_s^{**})> u_n(x)]-P[ M(X, K_s^{*})>u_n(\min\{x,y\})]-P[M(X, K_s^{**})> u_n(y)]\\
&&\ \ \ + P[M(\widehat{X}, K_s^{**})> u_n(x), M(X, K_s^{*})>u_n(\min\{x,y\})]\\
&&\ \ \ +P[M(\widehat{X}, K_s^{**})> u_n(x), M(X, K_s^{**})> u_n(y)]\\
&&\ \ \ +P[ M(X, K_s^{*})>u_n(\min\{x,y\}), M(X, K_s^{**})> u_n(y)]\\
&&\leq 1-\sharp(K_s^{**})(1-F(u_n(x)))-\sharp(K_s^{*})(1-F(u_n(\min\{x,y\})))-\sharp(K_s^{**})(1-F(u_n(y)))\\
&&\ \ \ +2t(1-F(u_n(\min\{x,y\})))^2+t\sum_{j=2}^{t} P(A_{s1},A_{sj})\\
&&= 1-t[1-F(u_n(\min\{x,y\}))]-(t-\sum_{j\in K_s}^{}\alpha_j )[1-F(u_n(\max\{x,y\}))]\\
&&\ \ \ +2t(1-F(u_n(\min\{x,y\})))^2+t\sum_{j=2}^{t} P(A_{s1},A_{sj})\\
&&= 1-t[1-F(u_n(\min\{x,y\}))]-(1-\frac{r}{2^k} )t[1-F(u_n(\max\{x,y\}))]\\
&&\ \ \ +\left  [  \frac{\sum_{j\in K_s}^{}\alpha_j  }{t} -\frac{r}{2^k}\right ]t[1-F(u_n(\max\{x,y\}))]\\
&&\ \ \ +2t(1-F(u_n(\min\{x,y\})))^2+t\sum_{j=2}^{t} P[A_{s1}, A_{sj}].
\end{eqnarray*}
Similarly, we have
\begin{eqnarray*}
&&P[M(X(\boldsymbol{\alpha}   ), K_s)\le u_n(x), M(X, K_s)\le u_n(y)]\\
&&\geq1- t[1-F(u_n(\min\{x,y\}))]+(t-\sum_{j\in K_s}^{}\alpha_j )[1-F(u_n(\max\{x,y\}))]\\
&&= 1-t[1-F(u_n(\min\{x,y\}))]-(1-\frac{r}{2^k} )t[1-F(u_n(\max\{x,y\}))]\\
&&\ \ \ +\left  [  \frac{\sum_{j\in K_s}^{}\alpha_j  }{t} -\frac{r}{2^k}\right ]t[1-F(u_n(\max\{x,y\}))].
\end{eqnarray*}
The proof of the lemma is complete.

\textbf{The proof of Theorem 2.1}. For random variable $\lambda$ such that $0 \leq \lambda \leq 1$ a.s., put
$$B_{r,l}=\begin{Bmatrix}
   w:\lambda (w)\in \left\{\begin{matrix}\left [0,\frac{1}{2^l}   \right ], &r=0
   \\ \left (\frac{r}{2^l} ,\frac{r+1}{2^l}\right],  &0<r<2^l-1
\end{matrix}\right.
\end{Bmatrix}$$
and
$$B_{r, l, \boldsymbol{\alpha}, n}=\{w: \varepsilon_j(w)= \alpha_j, 1 \leq j \leq n\}\cap B_{r, l}.$$
We split the proof into five steps. In the first step, we show that
\begin{eqnarray}
\label{Pr3.1.1}
&&\sum_{r=0}^{2^k-1} \sum_{\boldsymbol{\alpha}  \in \{0, 1\}}E\bigg|P[M_n(X(\boldsymbol{\alpha} ))\le u_n(x), M_n(X)\le u_n(y)]\nonumber\\
&&\ \ \ -\prod_{s=1}^{k}P[M(X(\boldsymbol{\alpha}   ), K_s)\le u_n(x), M(X, K_s)\le u_n(y)]\bigg|  I_{[B_{r, k, \boldsymbol{\alpha} , n}]}\rightarrow0,
\end{eqnarray}
as $n\rightarrow\infty$ and then $k\rightarrow\infty$. By Lemma 3.2, the sum in (\ref{Pr3.1.1}) is bounded above by
$$ 2(k-1) \alpha_{n,l} + (\frac{4kl}{n}+\frac{2}{k}) n [1-F(u_n(\min\{x,y\}))],$$
which tends to $0$ as $n\rightarrow\infty$ and then $k\rightarrow\infty$, since $\alpha_{n,l}\rightarrow0$ and
$n [1-F(u_n(\min\{x,y\}))]\rightarrow -\ln G(\min\{x,y\})$, as $n\rightarrow\infty$.

In the second step, we show
\begin{eqnarray}
\label{Pr3.1.2}
&&\sum_{r=0}^{2^k-1} \sum_{\boldsymbol{\alpha}  \in \{0, 1\}^n }E\bigg|\prod_{s=1}^{k}P[M(X(\boldsymbol{\alpha}), K_s)\le u_n(x), M(X, K_s)\le u_n(y)]\nonumber\\
&&\ \ \ -\left[1-\frac{n[1-F(u_n(\min\{x,y\}))]+(1-\frac{r}{2^k} )n[1-F(u_n(\max\{x,y\}))]}{k}\right ]^{k}  \bigg|  I_{[B_{r, k, \boldsymbol{\alpha} , n}]}\rightarrow0,\nonumber\\
\end{eqnarray}
as $n\rightarrow\infty$ and then $k\rightarrow\infty$.
Noting that
\begin{eqnarray}
 \label{Pr3.1.3}
\left |\prod_{s=1}^{k} a_s-\prod_{s=1}^{k} b_s\right |\le \sum_{s=1}^{k} \left |a_s-b_s  \right | ,
\end{eqnarray}
for all $a_s,b_s\in[0,1]$, by Lemma 3.3, we can bound the sum in (\ref{Pr3.1.2}) by
\begin{eqnarray}
\label{Pr3.1.4}
&&\Sigma:=\sum_{r=0}^{2^k-1} \sum_{s=1}^{k}E\left| \frac{\sum_{j\in K_s }^{}\varepsilon_j  }{t} -\frac{r}{2^k}\right|I_{[B_{r, k}]}\frac{n[1-F(u_n(\max\{x,y\}))]}{k}\nonumber \\
&&\ \ \ +2n(1-F(u_n(\min\{x,y\})))^2+n\sum_{j=2}^{t} P[A_{11},A_{1j}].
\end{eqnarray}
It follows from the assumption (\ref{2.1}) that
\begin{eqnarray*}
\frac{S_{st}}{st}\stackrel{p}\longrightarrow\lambda\ \mbox{and}\ \frac{S_{(s-1)t}}{(s-1)t}\stackrel{p}\longrightarrow\lambda,
\end{eqnarray*}
as $t\rightarrow\infty$,
which combined with the dominated convergence theorem yields
$$E\bigg | \frac{S_{st}}{st} -\lambda \bigg|\rightarrow0\ \ \mbox{and}\ \ E \bigg |\frac{S_{(s-1)t}}{(s-1)t} -\lambda  \bigg |\rightarrow0,$$
as $t\rightarrow\infty$.
Thus we have
\begin{eqnarray}
 \label{Pr3.1.5}
 &&\sum_{r=0}^{2^k-1}E\bigg| \frac{\sum_{j\in K_s }^{}\varepsilon_j  }{t} -\frac{r}{2^k}\bigg|I_{[B_{r, k}]}\nonumber \\
 &&\leq E\bigg| \frac{\sum_{j\in K_s }^{}\varepsilon_j  }{t} -\lambda\bigg|+\sum_{r=0}^{2^k-1}E\bigg|\lambda-\frac{r}{2^k}\bigg|I_{[B_{r, k}]}\nonumber \\
 &&\leq E\bigg |\frac{S_{st}-S_{(s-1)t}}{t}-\lambda \bigg | +\frac{1}{2^k}\nonumber \\
 &&\leq  E\bigg | s \left(\frac{S_{st}}{st} -\lambda \right)+ (s-1) \left(\frac{S_{(s-1)t}}{(s-1)t} -\lambda \right) \bigg | +\frac{1}{2^k}\nonumber \\
 &&\leq s E\bigg | \frac{S_{st}}{st} -\lambda \bigg|+ (s-1)E \bigg |\frac{S_{(s-1)t}}{(s-1)t} -\lambda  \bigg | +\frac{1}{2^k}\nonumber  \\
 &&\leq o(1)+ \frac{1}{2^k},
\end{eqnarray}
as $t\rightarrow\infty$, which combinied with condition $D'(u_{n})$ and that facts that $n [1-F(u_n(\max\{x,y\}))]\rightarrow -\ln G(\max\{x,y\})$, as $n\rightarrow\infty$, implies
\begin{eqnarray*}
\Sigma\leq \frac{-\ln G(\max\{x,y\})}{2^k},
\end{eqnarray*}
as $n\rightarrow\infty$. Letting $k\rightarrow\infty$, we finished the proof of (\ref{Pr3.1.2}).

In the third step, we prove
\begin{eqnarray}
 \label{Pr3.1.6}
&&\sum_{r=0}^{2^k-1}\sum_{\boldsymbol{\alpha}\in\{0,1\}^n}E\bigg| \left[1- \frac{n[1-F(u_n(\min\{x,y\}))]+(1-\frac{r}{2^k} )n[1-F(u_n(\max\{x,y\}))]}{k}   \right ]^{k} \nonumber  \\
 &&\ \ \ \ \ \ \ \   -\left[1- \frac{n[1-F(u_n(\min\{x,y\}))]+(1-\lambda )n[1-F(u_n(\max\{x,y\}))]}{k}   \right ]^{k}  \bigg|I_{[B_{r, k, \boldsymbol{\alpha},n}]} \rightarrow0,\nonumber\\
\end{eqnarray}
as $n\rightarrow\infty$ and then $k\rightarrow\infty$.
Using (\ref{Pr3.1.3}) again, we can see the left hand side of (\ref{Pr3.1.6}) is bounded above by
\begin{eqnarray*}
&& \sum_{r=0}^{2^k-1} E\left |\frac{r}{2^k} -\lambda   \right | I_{[B_{r, k}]}n[1-F(u_n(\max\{x,y\}))]
\leq \frac{n[1-F(u_n(\max\{x,y\}))]}{2^k},
\end{eqnarray*}
which tends to $0$ as $n\rightarrow\infty$ and then $k\rightarrow\infty$, since $n[1-F(u_n(\max\{x,y\}))]\rightarrow -\ln G(\max\{x,y\})$, as $n\rightarrow\infty$. This completes the proof of (\ref{Pr3.1.6}).

In the fourth step, we show
\begin{eqnarray}
 \label{Pr3.1.7}
&&\sum_{r=0}^{2^k-1}\sum_{\boldsymbol{\alpha}\in\{0,1\}^n}E\bigg|\left[1- \frac{n[1-F(u_n(\min\{x,y\}))]+(1-\lambda )n[1-F(u_n(\max\{x,y\}))]}{k} \right ]^{k}\nonumber  \\
 &&\ \ \ \ \ \ \ \   - \left [1-\frac{-\ln G(\min\{x, y\})-\ln G^{1-\lambda }(\max\{x, y\})}{k}   \right ]^k\bigg|I_{[B_{r, k, \boldsymbol{\alpha},n}]}\rightarrow0,
\end{eqnarray}
as $n\rightarrow\infty$ and then $k\rightarrow\infty$. This result follows easily from  (\ref{Pr3.1.3}) and  the facts that
$n[1-F(u_n(\max\{x,y\}))]\rightarrow -\ln G(\max\{x,y\})$ and $n[1-F(u_n(\min\{x,y\}))]\rightarrow -\ln G(\min\{x,y\})$
as $n\rightarrow\infty$.

In the last step, letting $k\rightarrow\infty$,  we have
$$\bigg|E{\left [1-\frac{-\ln G(\min\{x, y\})-\ln G^{1-\lambda }(\max\{x, y\})}{k}   \right ]}^k -G(\min\{x, y\})E G^{1-\lambda }(\max\{x, y\})\bigg|\rightarrow0.$$
Thus, combining the above five steps, we get the desired result
$$\lim_{n \to \infty} P[M_n(X(\boldsymbol{\varepsilon} ))\le u_n(x), M_n(X)\le u_n(y)]=G(\min\{x, y\})EG^{1-\lambda }(\max\{x, y\}).$$

\textbf{The proof of Theorem 2.2}. For stationary Gaussian sequences with covariance function $r_{n}$ satisfying $r_{n}\log n\rightarrow0$ as $n\rightarrow\infty$, conditions $D(u_{n},v_{n})$ and  $D'(u_{n})$ hold, see Lemma 4.4.1 in Leadbetter et al. (1983) for details.  Hence,
 Theorem 2.2 follows from Theorem 2.1.

\textbf{The proof of Theorem 2.3}. By the similar arguments as in the proof of Lemma 4.4.1 in Leadbetter et al. (1983) and using a comparison lemma for $chi$-random variables (see e.g., Lemma 3.2 in Song et al. (2022)) to replace the Normal Comparison Lemma, we can prove that conditions $D(u_{n},v_{n})$ and  $D'(u_{n})$ hold for the $chi$-random variables defined in Theorem 2.3. Hence, Theorem 2.3 also follows from Theorem 2.1.

\textbf{The proof of Theorem 2.4}. By the similar arguments as in the proof of Lemma 4.4.1 in Leadbetter et al. (1983) and using a comparison lemma for Gaussian order statistics (see e.g., Theorem 2.4 in  D\c{e}bicki et al. (2017)) to replace the Normal Comparison Lemma, we can prove that conditions $D(u_{n},v_{n})$ and  $D'(u_{n})$ hold also for Gaussian order statistics random variable defined in Theorem 2.4. Hence, Theorem 2.4 follows again from Theorem 2.1.

\bigskip


\end{document}